\documentclass[12pt,reqno]{amsart}
\usepackage[margin=1in]{geometry}
\usepackage[utf8]{inputenc}
\usepackage{setspace}
\usepackage{graphicx}
\usepackage{xcolor}
\usepackage{amsmath, mathtools}
\usepackage{amsthm}
\usepackage{enumitem}
\usepackage{caption}
\usepackage{subcaption}

\usepackage{todonotes}
\usepackage{tikz}
\usepackage[normalem]{ulem}

\newtheorem{theorem}{Theorem}

\newtheorem{observation}{Observation}

\newtheorem{definition}{Definition}
\newtheorem{lemma}{Lemma}

\newtheorem{claim}{Claim}

\newtheorem{conjecture}{Conjecture}

\usepackage{url}
\usepackage{hyperref}
\hypersetup{
    colorlinks,
    linkcolor={red!60!black},
    citecolor={green!60!black},
    urlcolor={blue!60!black}
}

\tikzstyle{vert}=[shape=circle,draw=black,fill=black, inner sep=.75mm]
\tikzstyle{fixed}=[shape=rectangle,draw=black,fill=white, inner sep=1.2mm]
\tikzstyle{uncolored}=[dashed,thick]
\tikzstyle{uncolored2}=[dotted,thick]
\tikzstyle{red?}=[dashed,thick,color=red]
\tikzstyle{blue?}=[dashed,thick,color=blue]
\tikzstyle{purple?}=[dashed,thick,color=purple]
\tikzstyle{purple}=[solid,thick,color=purple]
\tikzstyle{red}=[solid,thick,color=red]
\tikzstyle{blue}=[solid,thick,color=blue]
\tikzstyle{green}=[solid,thick,color=green]
\tikzstyle{green?}=[dashed,thick,color=green]

\newcommand{\rbrac}[1]{\left(#1\right)} 

\newcommand{\mc}[1]{\mathcal{#1}}
\newcommand{\tbf}[1]{\emph{#1}}

\newcommand{\eps}{\varepsilon}

\newcommand{\nn}{\nonumber}

\def\sm{\setminus}

\allowdisplaybreaks

\title{Generalized Ramsey numbers at the linear and quadratic thresholds}

\author{Patrick Bennett}
\address{Department of Mathematics, Western Michigan University, Kalamazoo, MI, USA}
\thanks{The first author was supported in part by Simons Foundation Grant \#426894.}
\email{\tt patrick.bennett@wmich.edu}

\author{Ryan Cushman}\thanks{}
\address{Department of Mathematics,
University of Wisconsin-Eau Claire, Eau Claire, WI, USA} 
\email{\tt cushmarj@uwec.edu}

\author{Andrzej Dudek}
\address{Department of Mathematics, Western Michigan University, Kalamazoo, MI, USA}
\thanks{The third author was supported in part by Simons Foundation Grant \#522400.}
\email{\tt andrzej.dudek@wmich.edu}

\begin{document}

\maketitle 

\begin{abstract}
     The generalized Ramsey number $f(n, p, q)$ is the smallest number of colors needed to color the edges of the complete graph $K_n$ so that every $p$-clique spans at least $q$ colors. Erd\H{o}s and Gy\'arf\'as showed that $f(n, p, q)$ grows linearly in $n$ when $p$ is fixed and $q=q_{\text{lin}}(p):=\binom p2-p+3$. Similarly they showed that $f(n, p, q)$ is quadratic in $n$ when $p$ is fixed and $q=q_{\text{quad}}(p):=\binom p2-\frac p2+2$. In this note we improve on the known estimates for $f(n, p, q_{\text{lin}})$ and $f(n, p, q_{\text{quad}})$. Our proofs involve establishing a significant strengthening of a previously known connection between $f(n, p, q)$ and another extremal problem first studied by Brown, Erd\H{os} and S\'os, as well as building on some recent progress on this extremal problem by Delcourt and Postle and by Shangguan. Also, our upper bound on $f(n, p, q_{\text{lin}})$ follows from an application of the recent forbidden submatchings method of Delcourt and Postle. 
\end{abstract}

\section{Introduction}

Erd\H{o}s and Shelah~\cite{E75} first considered the following generalization of classical Ramsey problem. 
\begin{definition}
Fix  integers $p, q$ such that $p \ge 3$ and $2 \le q \le \binom p2$. A \emph{$(p, q)$-coloring} of $K_n$ is a coloring of the edges of $K_n$ such that every $p$-clique has at least $q$ distinct colors among its edges. The generalized Ramsey number $f(n, p, q)$ is the minimum number of colors such that $K_n$ has a $(p, q)$-coloring.
\end{definition}

Erd\H{o}s and Gy\'arf\'as \cite{EG97} systematically studied $f(n, p, q)$ for fixed $p, q$ as $n \rightarrow \infty$. In this paper all asymptotic statements are as $n \rightarrow \infty$. Among other results, Erd\H{o}s and Gy\'arf\'as \cite{EG97} proved that for arbitrary $p$ and 
\[
q=q_{\text{lin}}(p):=\binom p2-p+3,
\]
$f(n,p,q)$ is linear, but $f(n,p,q-1)$ is sublinear. Similarly, they showed in~\cite{EG97} that for 
\[
q=q_{\text{quad}}(p):=\binom p2-\lfloor p/2\rfloor+2,
\]
$f(n,p,q)$ is quadratic, but $f(n,p,q-1)$ is subquadratic. Thus for fixed $p$, we call the value $q_{\text{lin}}$ the \textit{linear threshold} and $q_{\text{quad}}$ the \textit{quadratic threshold}. The main goal of this note is to estimate $f(n, p, q)$ when $q$ is at the linear or quadratic threshold. In terms of explicit general bounds, we prove the following. 

\begin{theorem}\label{thm:est}
For all $p \ge 3$ we have
    \begin{equation}\label{eqn:estlin}
    \frac{3p-7}{4p-10}n +o(n)  \le   f(n, p, q_{\text{lin}}) \le n+o(n).   
    \end{equation}
For even $p \ge 6$ we have
    \begin{equation}\label{eqn:estquad}
        \frac{2p-7}{5p-18}n^2 +o(n^2)  \le   f(n, p, q_{\text{quad}}) \le \frac 5{12} n^2+o(n^2).
    \end{equation}
\end{theorem}
\noindent
Some of the content of Theorem \ref{thm:est} was also recently obtained independently by G\'omez-Leos, Heath, Parker, Schwieder, and Zerbib \cite{GHPSZ}.

Since the initial investigation by Erd\H{o}s and Gy\'arf\'as \cite{EG97}, the asymptotic behavior of $f(n,p,q)$ has attracted a considerable amount of attention. See, for example, the recent paper of the first author, third author and English \cite{BDE22} for some history of the problem. 
However, except for the trivial case of $f(n, 3, 3) = n+O(1)$, there have only been two results where $f(n,p,q)$ is known with a $(1+o(1))$ multiplicative error. Erd\H{o}s and Gy\'arf\'as \cite{EG97} stated that it ``can be easily determined''  that
\begin{equation}\label{eqn:6,14}
    f(n, 6, 14)=\frac{5}{12}  n^2+O(n).
\end{equation}
 More recently, the present authors with Pra\l at~\cite{BCDP22}  proved that $f(n, 4, 5) = \frac 56 n + o(n)$. In this note we provide a proof for \eqref{eqn:6,14} and also obtain $f(n, 6, 14)$ exactly when $n \equiv 1, 4 \pmod{12}$ (see Theorem \ref{thm:6,14}). We also obtain two more more explicit and asymptotically sharp estimates for generalized Ramsey numbers at the quadratic threshold. 
\begin{theorem}\label{thm:8,26}
We have
    \[
    f(n, 8, 26) = \frac{9}{22} n^2 +o(n^2) \quad\text{and}\quad f(n, 10, 42) = \frac{5}{12}n^2+o(n^2).
    \]
\end{theorem}

The proofs of Theorems \ref{thm:est} and \ref{thm:8,26} will involve establishing certain connections between $f(n, p, q)$ and the following extremal problem first studied by Brown, Erd\H{os} and S\'os~\cite{BES1973}. 
\begin{definition}
 Let $\mc{H}$ be an $r$-uniform hypergraph. A \emph{$(s, k)$-configuration} in $\mc{H}$ is a set of $s$ vertices inducing $k$ or more edges. We say $\mc{H}$ is \emph{$(s, k)$-free} if it has no $(s, k)$-configuration. Let $F^{(r)}(n; s, k)$ be the largest possible number of edges in an $(s, k)$-free $r$-uniform hypergraph with $n$ vertices. In terms of classical extremal numbers, 
\[
  F^{(r)}(n; s, k) =  \mbox{ex}_{r}(n, \mc{G}_{s, k}),
\]
 where $\mc{G}_{s, k}$ is the family of all $r$-uniform hypergraphs on $s$ vertices and $k$ edges.
\end{definition}
\noindent 
In fact, three of the four explicit bounds in Theorem \ref{thm:est} follow by first bounding $f(n, p, q)$ implicitly in terms of some values $F^{(r)}(n; s, k)$ and then using explicit bounds on the latter. Thus, further improvements on the estimates for $F^{(r)}(n; s, k)$ would in some cases automatically give improved estimates for $f(n, p, q)$. In the case of the quadratic threshold (and even $p$), we actually show that the problem of asymptotically estimating $f(n, p, q)$ completely reduces to asymptotically estimating a certain value $F^{(r)}(n; s, k)$.

\begin{theorem}\label{thm:quadexact}
For all even $p \ge 6$ we have
\begin{equation}\label{eqn:quadexact}
    \lim_{n \rightarrow \infty} \frac{f\rbrac{n, p, q_{\text{quad}}}}{n^2} = \frac12 - \lim_{n \rightarrow \infty} \frac{F^{(4)}\rbrac{n; p, \frac p2-1}}{n^2}.
\end{equation}
    In particular, the limit on the left exists. Furthermore, there exist asymptotically optimal $\rbrac{p, q_{\text{quad}}}$-colorings that use no color more than twice. 
\end{theorem}
\noindent
It is perhaps surprising that we need not use any color more than twice. Indeed a $(p, q_{\text{quad}})$-coloring is allowed to use a color up to $\frac p2 -1$ times, and it would seem more efficient to use the same color as many times as possible. 

The existence of the limit on the right-hand side of \eqref{eqn:quadexact} was proved only recently by Shangguan \cite{S22}. Shangguan's proof generalizes another recent result by Delcourt and Postle~\cite{DP22}, which resolved a conjecture from Brown, Erd\H{os} and S\'os \cite{BES1973} regarding the existence of a similar limit involving the function $F^{(3)}$ for 3-uniform hypergraphs. In particular, Delcourt and Postle \cite{DP22} proved the existence (for fixed $\ell \ge 3$) of the limit
\[
    \lim_{n \rightarrow \infty} \frac{F^{(3)}\rbrac{n; \ell, \ell-2}}{n^2}.
\]
\noindent
Interestingly, the proofs of Delcourt and Postle \cite{DP22} and Shangguan \cite{S22} do not seem to shed much light (at least, not as much as one might hope) on how to actually find the limits whose existence they establish. However, these limits are known in two cases relevant to us. In particular, it is known due to Shangguan and Tamo \cite{ST2019} that
\begin{equation}\label{eqn:chong}
    F^{(4)}(n; 8, 3) = \frac{1}{11} n^2 +o(n^2).
\end{equation}
It is also known due to Glock, Joos, Kuhn, Kim, Lichev and Pikhurko \cite{GJKKLP2022} that
\begin{equation}\label{eqn:oleg}
    F^{(4)}(n; 10, 4) = \frac{1}{12}n^2 +o(n^2).
\end{equation}
Thus, Theorem \ref{thm:8,26} follows from Theorem \ref{thm:quadexact} together with \eqref{eqn:chong} and \eqref{eqn:oleg}.

The lower bound on $f(n, p, q_{\text{lin}})$ in Theorem \ref{thm:est} is similar to the quadratic case, in the sense that it follows from an upper bound on $F^{(3)}(n;p, p-2)$. In particular we prove that 
\begin{theorem}\label{thm:linlower}
For all $p\ge 3$ we have
    \begin{equation}\label{eqn:linlower}
    \liminf_{n \rightarrow \infty} \frac{f\rbrac{n, p, q_{\text{lin}}}}{n} \ge 1 - \lim_{n \rightarrow \infty} \frac{F^{(3)}\rbrac{n; p,  p-2}}{n^2}.
\end{equation}
\end{theorem}
\noindent
In light of Theorem \ref{thm:quadexact}, one might suspect that there is a matching upper bound for \eqref{eqn:linlower}, but unfortunately this is not the case. Indeed, Glock \cite{G19} proved that 
\[
\lim_{n \rightarrow \infty} \frac{F^{(3)}\rbrac{n; 5,  3}}{n^2} = \frac 15,
\]
which together with \eqref{eqn:linlower} yields $f(n, 5, 8) \ge \frac 45 n + o(n)$. However this lower bound is not close to the truth. Indeed, recently it was proved that $f(n, 5, 8)= \frac 67 n + o(n)$. The lower bound was proved by G\'omez-Leos, Heath, Parker, Schwieder, and Zerbib \cite{GHPSZ} and the upper bound by the current authors \cite{BCD24}.

\subsection{Comparison to previous bounds}

Here we compare the bounds in Theorem \ref{thm:est} to what was previously known. For the linear threshold, Erd\H{o}s and Gy\'arf\'as~\cite{EG97} showed that 
\begin{equation}\label{eqn:linprev}
   (n-1)/(p-2) \le f(n, p, q_{\text{lin}}) \le c_p n 
\end{equation}
for some coefficient $c_p$. The lower bound in \eqref{eqn:linprev} follows from the simple fact that in a $(p, q_{\text{lin}})$-coloring each vertex is adjacent to at most $p-2$ edges of each color. The upper bound in \eqref{eqn:linprev} follows from the Local Lemma. The constant $c_p$ is not explicitly discussed in~\cite{EG97} but it is easy to see from their proof that $c_p \rightarrow \infty$ as $p \rightarrow \infty$. 
Thus we see that in Theorem \ref{thm:est}, \eqref{eqn:estlin} is a  significant improvement on previous bounds. Indeed, the gap between the coefficients in \eqref{eqn:linprev} grows without bound with $p$, whereas the coefficients in \eqref{eqn:estlin} are always between $3/4$ and $1$.

Likewise, for the quadratic threshold (and even $p$) the trivial bounds are
\[
\frac{\binom n2}{\frac p2 -1} \le f(n, p, q_{\text{quad}}) \le \binom n2.
\]
The upper bound follows since we can give every edge its own color, and the lower bound follows from the fact that each color must be used at most $\frac p2 -1$ times. Thus we see that \eqref{eqn:estquad} in Theorem \ref{thm:est} is a significant improvement. 

\subsection{Structure of the note}

The structure of this note is as follows. In Section~\ref{sec:quad} we address the quadratic threshold. We start with a proof of a more precise version of \eqref{eqn:6,14}. We go on to prove Theorem \ref{thm:quadexact} and \eqref{eqn:estquad} from Theorem \ref{thm:est}. In Section~\ref{sec:lin} we address the linear threshold. There we prove Theorem \ref{thm:linlower} and \eqref{eqn:estlin} from Theorem \ref{thm:est}.

\section{Quadratic Threshold}\label{sec:quad}

In this section we address the quadratic threshold. First we introduce some terminology. Suppose we are given a coloring of the edges of $K_n$. For a set of vertices $S$, let $c(S)$ be the number of colors appearing on edges within $S$, and let $r(S)$ be $\binom{|S|}{2}-c(S)$. We call $r(S)$ the \tbf{number of color repetitions (or just repeats) in $S$}. Sometimes it may help the reader to imagine counting $r(S)$ by examining each edge of $S$ in some order and counting a repeat whenever we see a color we have already seen.

\subsection{Estimating $f(n, 6, 14)$} In this subsection we state and prove our more precise result for $f(n, 6, 14)$. As we noted, Erd\H{o}s and Gy\'arf\'as \cite{EG97} stated that $f(n, 6, 14)=\frac{5}{12}  n^2+O(n)$ without proof. To help the reader gain familiarity with the concepts in this paper, we present a proof of a more precise version of this result.

\begin{theorem}\label{thm:6,14}
We have
\[
\frac{5}{6}\binom{n}{2} \le f(n, 6, 14) \le \frac{5}{6}\binom{n}{2}+O(n).
\]
Furthermore, the lower bound above is the exact value of $f(n, 6, 14)$ whenever $n$ is congruent to $1$ or $4$ modulo $12$. 
\end{theorem}

\begin{proof}
Starting with the lower bound, suppose we have any $(6,14)$-coloring. Since $\binom{6}{2}=15$,
any set of $6$ vertices is allowed to have only one repeat, which implies that we cannot have $3$ edges of the same color. Indeed, taking the union of these edges would be a set of at most $6$ vertices with more than one repeated color.  
This also means that there can be at most one monochromatic path on three vertices $P_3$, since the union of two of them would be a set of at most 6 vertices with at least two repeats. If our coloring contains a monochromatic $P_3$, then we remove it and get a coloring of $K_{n-3}$. So we have a $(6, 14)$-coloring of $K_{n'}$ with $n' \in \{n-3, n\}$ with no monochromatic $P_3$.

Suppose the color $c$ is used twice, say on the (nonincident) edges $ab$ and $xy$. Then the other four edges in $\{a, b, x, y\}$ must all have different colors which are only used once in the whole graph. Let $C_1$ be the set of colors used once and $C_2$ the colors used twice. For each $c \in C_2$ let $K_c$ be the set of 4 vertices consisting of both endpoints of both edges of color $c$. Note that for $c, c' \in C_2$ we have $|K_c \cap K_{c'}| \le 1$, since otherwise $K_c \cup K_{c'}$ is a set of at most 6 vertices with too many repeats.  Thus the sets $K_c$ induce edge-disjoint 4-cliques. Thus, if we did not remove any $P_3$, we have that $|C_2|$, the number of such cliques, is at most 
\begin{equation}\label{eq:f614-lb-1}
\frac16\binom{n}2.
\end{equation}
\noindent
On the other hand, if we did remove a $P_3$, this would contribute one additional color to $C_2$ along with the restriction on the $K_c$. From our discussion above, we note that this $P_3$ is vertex disjoint from all the $K_c$ and does not share a color with any other edges. Thus, in this case, $|C_2|$ is at most 
\begin{equation}\label{eq:f614-lb-2}
1 + \frac16\binom{n-3}{2}.
\end{equation}
\noindent
But since \eqref{eq:f614-lb-1} is at least \eqref{eq:f614-lb-2} for $n\ge 4$, we conclude that the number of colors used is at least
$$
|C_1| + |C_2| = \left(\binom{n}{2}-2|C_2|\right) + |C_2|=\binom{n}{2}-|C_2| \ge \frac56 \binom{n}2.
$$
\noindent
Thus we are done with the lower bound for Theorem \ref{thm:6,14}. We move on to the upper bound.

 If $n \equiv 1$ or $4 \mod 12$, then we are guaranteed a perfect packing of $\frac16\binom{n}{2}$ edge-disjoint 4-cliques by Hanani's result~\cite{H61}.
Then for each clique in the packing, color two nonadjacent edges the same color and give a unique color to the remaining edges. Since we use exactly 5 colors for each clique, we use exactly $\frac56\binom{n}{2}$ colors to color all the edges. 

Otherwise, let $i = (n \mod 12)$ and partition the vertices into $K_{n-i+1} \cup K_{i-1}$, and find a perfect packing of edge-disjoint 4-cliques for $K_{n-i+1}$. Follow the same coloring as above for the perfect packing, and then color the remaining $(n-i+1)(i-1)+\binom{i-1}{2} = O(n)$ edges with a different color for each edge.  
Thus, we use $\frac56\binom{n-i+1}{2}+O(n)=\frac56\binom{n}{2} + O(n)$ colors.

Notice that in either case, the resulting coloring satisfies the $(6,14)$-coloring condition. If not, then there exists a set $S$ of $6$ vertices with more than 2 repeated colors. In our coloring, this means that $S$ must contain two cliques from the packing. But since the cliques must be edge-disjoint, this implies that $|S|\ge 7$, a contradiction. 
\end{proof}

\subsection{Proof of Theorem \ref{thm:quadexact}}
In this subsection we will prove Theorem \ref{thm:quadexact} after some discussion. We consider the case of $(p, q)$-coloring, where 
\[
p=2\ell \quad{\text{and}}\quad q= q_{\text{quad}}(p) = \binom{2\ell}{2}-\ell+2. 
\]
This choice of parameters allows using a color $\ell-1$ times but not $\ell$ times. Erd\H{o}s and Gy\'arf\'as \cite{EG97} showed that for this choice of parameters $f(n, p, q)$ is quadratic in $n$. Of course the upper bound $f(n, p, q) \le \binom n2$ is trivial, but \cite{EG97} also gives a nontrivial upper bound of $(1/2-\eps)n^2$ for some $\eps  >0$. Specifically, Erd\H{o}s and Gy\'arf\'as \cite{EG97} used a 4-uniform $(2 \ell, \ell-1)$-free hypergraph $\mc{H}$ to give an appropriate coloring. Crucially, every color repetition in the coloring corresponds to a hyperedge of $\mc{H}$. Specifically each color is used at most twice, and for any color used on two edges, the union of those two edges is a hyperedge of $\mc{H}$. The existence of a suitable hypergraph had already been established by Brown, Erd\H{o}s and S\'os~\cite{BES1973}. The same basic connection between $(p, q)$-coloring near the quadratic threshold and 4-uniform $(s, k)$-free hypergraphs (for the appropriate $s, k$) was exploited by S\'{a}rk\"{o}zy and Selkow \cite{SS03} and again by Conlon, Gishboliner, Levanzov and Shapira \cite{CGLS23}. However, this connection as it was used in \cite{CFLS15,EG97,SS03} is not precise enough to prove Theorem \ref{thm:quadexact}. Indeed, all these previous results give away a constant factor in the main term of their estimate of $f(n, p, q)$, while we want an asymptotically tight estimate. Thus, we will have to significantly refine these previously established connections between the Erd\H{o}s-Gy\'arf\'as coloring problem and the Brown-Erd\H{o}s-S\'os packing problem.

Now we will define some functions related to $F^{(4)}(n; 2\ell, \ell-1)$. The first one relaxes the problem to multi-hypergraphs.

\begin{definition}
 Let $\mc{H}$ be an $r$-uniform \emph{multi-hypergraph}, meaning that $\mc{H}$ can have edges with multiplicity (but each edge has $r$ distinct vertices). A \emph{$(s, k)$-configuration} in $\mc{H}$ is a set of $s$ vertices inducing $k$ or more edges (counted with multiplicity). Let $G^{(r)}(n; s, k)$ be the largest possible number of edges in an $(s, k)$-free $r$-uniform multi-hypergraph with $n$ vertices.
\end{definition}
\noindent
Next we define a function that restricts the extremal problem for $F^{(4)}(n; 2\ell, \ell-1)$ to a smaller family of hypergraphs.

\begin{definition}\label{def:H}
  Let $H^{(4)}(n; 2\ell, \ell-1)$ be the largest possible number of edges in a $4$-uniform hypergraph $\mc{H}$ on $n$ vertices which satisfies the following conditions:
 \begin{enumerate}
     \item \label{item:1} $\mc{H}$ is $(2\ell, \ell-1)$-free,
     \item \label{item:2} $\mc{H}$ is $(2i+1, i)$-free for $i=2, \ldots, \ell-2$, and 
     \item \label{item:3} for every vertex $v$ of $\mc{H}$, either $v$ has degree $0$ or degree at least $\ell-1$. 
 \end{enumerate}
\end{definition}
\noindent
Using Shangguan's notation \cite{S22}, our function $H^{(4)}(n; 2\ell, \ell-1)$ defined above is the same as what Shangguan refers to as $f_r^{(t)}(n; er-(e-1)k, e)$, where $4$ is substituted for $r$, $2$ for $k$, $2$ for $t$, and $\ell-1$ for $e$. Since $H^{(4)}$ is a restriction and $G^{(4)}$ is a relaxation, we have
\[
    H^{(4)}(n; 2\ell, \ell-1) \le F^{(4)}(n; 2\ell, \ell-1)\le G^{(4)}(n; 2\ell, \ell-1).
\]
Shangguan \cite{S22} proved (see Lemma 5.5 and the discussion above it) that 
\begin{lemma}[Lemma 5.5 in \cite{S22}]\label{lem:shangguan}
\begin{equation}\label{eqn:FH}
   \lim_{n \rightarrow \infty} \frac{H^{(4)}(n; 2\ell, \ell-1)}{n^2} = \lim_{n \rightarrow \infty} \frac{F^{(4)}(n; 2\ell, \ell-1)}{n^2}. 
\end{equation}
\end{lemma}
\noindent
Now we will easily see that $G$ is likewise asymptotically the same as the others. 
\begin{claim}\label{clm:FG}
\begin{equation}\label{eqn:FG}
    \lim_{n \rightarrow \infty} \frac{G^{(4)}(n; 2\ell, \ell-1)}{n^2} = \lim_{n \rightarrow \infty} \frac{F^{(4)}(n; 2\ell, \ell-1)}{n^2}.
\end{equation}
\end{claim}
\begin{proof}
    Let $\mc{H}$ be an extremal multi-hypergraph for the $G^{(4)}(n; 2\ell, \ell-1)$ problem, i.e., $\mc{H}$ has $G^{(4)}(n; 2\ell, \ell-1)$ edges and is $(2\ell, \ell-1)$-free. We form a new hypergraph $\mc{H}'$ by simply deleting all multiple edges in $\mc{H}$. Clearly  $\mc{H}'$ is $(2\ell, \ell-1)$-free, so it has at most $F^{(4)}(n; 2\ell, \ell-1)$ edges. 

    We show that $\mc{H}'$ has almost the same number of edges as $\mc{H}$. Indeed, suppose we enumerate all the multiple edges of $\mc{H}$, say $\{e_1, \ldots, e_a\}$ where $e_i$ has multiplicity $m_i \ge 2$. It is easy to see that $a \le \ell$ and each $m_i \le \ell$ (otherwise there is a $(2\ell, \ell-1)$-configuration). Thus, we remove at most $\ell^2$ edges from $\mc{H}$ to get $\mc{H}'$. Consequently, we have
    \[
  F^{(4)}(n; 2\ell, \ell-1)  \le G^{(4)}(n; 2\ell, \ell-1)\le F^{(4)}(n; 2\ell, \ell-1) + \ell^2
    \]
    and \eqref{eqn:FG} follows (recall we already knew that the limit on the right exists due to Lemma~\ref{lem:shangguan}). 
\end{proof}

Now we start to attack the lower bound for the coloring problem. 

\begin{claim}\label{clm:quadlower}
\[
    f\rbrac{n, p, q_{\text{quad}}} \ge \binom n2 - G^{(4)}(n; 2\ell, \ell-1).
\]
\end{claim}
We will prove this claim directly, by using a $\rbrac{ p, q_{\text{quad}}}$-coloring to construct a $(2\ell, \ell-1)$-free multi-hypergraph. Towards that end we define the following. 
\begin{definition}
    Consider any coloring $C$ of the edges of $K_n$. We say a 4-uniform hypergraph $\mc{H}$ is a \emph{repeat multi-hypergraph} for the coloring $C$ if it is formed as follows. $\mc{H}$ has the same vertex set as $K_n$. For each color $c$ used in the coloring, let $E(c)\neq \emptyset$ be the set of edges of color $c$ and let $e_c$ be some particular (arbitrary) edge of color $c$. Then $\mc{H}$ will have all the edges $\{e \cup e_c: e \in E(c) \sm \{e_c\} \}$. Of course, $e \cup e_c$ might only have 3 vertices (when we claimed $\mc{H}$ would be 4-uniform) but we fix this by arbitrarily adding vertices to edges of size 3. 
\end{definition}
Note that $\mc{H}$ can have multiple edges since a single set of 4 vertices can contain, say, two red edges and also two blue edges. Also, since the construction of $\mc{H}$ potentially involves some arbitrary choices (in particular, the choice of the edges $e_c$ as well as the choice of vertices used to enlarge 3-edges), in general a coloring $C$ may give rise to several possible repeat multi-hypergraphs $\mc{H}$. 

We now make the key observation about repeat multi-hypergraphs. Essentially it says that edges in $\mc{H}$ count color repetitions of $C$ ``faithfully,'' i.e.,  without under- or over-counting. 

\begin{observation}\label{obs:faithful}
    Let $\mc{H}$ be a repeat multi-hypergraph for the coloring $C$. Then for all sets of vertices $S$ we have
    \[
    r(S) = |E(\mc{H}[S])|.
    \]
\end{observation}

\begin{proof}[Proof of Observation \ref{obs:faithful}]
  Recall that each hyperedge of $\mc{H}$ contains $e \cup e_c$ for some color $c$ and some edge $e$ that has color $c$. Now if a set of vertices $S$ vertices spans $b$ hyperedges all corresponding to the same color $c$, then $S$ contains $e_i \cup e_c$ for $1 \le i \le b$ and some edges $e_1, \ldots, e_b$ which all have color $c$. In particular $S$ contains $b+1$ edges, namely $e_c, e_1, \ldots e_b$, which all have color $c$, i.e.,  $S$ spans $b$ repeats in the color $c$. Now if $S$ spans $b$ hyperedges (which now need not all correspond to the same color), we likewise conclude that $S$ spans $b$ repeats by simply summing over the colors. 
  \end{proof}




We are now ready to prove Claim \ref{clm:quadlower}.
\begin{proof}[Proof of Claim \ref{clm:quadlower}]
    Consider a $(p, q_{\text{quad}})$-coloring $C$ of $K_n$ that is optimal, i.e.,  uses $f\rbrac{n, p, q_{\text{quad}}}$ colors. In such a coloring, any set of $p=2\ell$ vertices spans at most $\ell-2$ repeats. Let $\mc{H}$ be a repeat multi-hypergraph for $C$. By Observation \ref{obs:faithful}, a $(2\ell, \ell-1)$-configuration in $\mc{H}$ would be a set of $2 \ell$ vertices spanning at least $\ell-1$ repeats. Since $C$ is a $(2\ell, \binom{2\ell}{2}-\ell+2)$-coloring, $\mc{H}$ is $(2\ell, \ell-1)$-free. In particular, $|E(\mc{H})| \le G^{(4)}(n; 2\ell, \ell-1).$ But now applying Observation \ref{obs:faithful} where $S$ is the set of all vertices we have $|E(\mc{H})| = \binom n2 - f\rbrac{n, p, q_{\text{quad}}}$ since $C$ uses $f\rbrac{n, p, q_{\text{quad}}}$ colors. The claim now follows from
    \[
    \binom n2 - f\rbrac{n, p, q_{\text{quad}}} = |E(\mc{H})| \le G^{(4)}(n; 2\ell, \ell-1)
    \]
\end{proof}

Next we attack the upper bound for the coloring problem. To get a bound that comes close to matching Claim \ref{clm:quadlower}, we will have to ``reverse'' the procedure we used to turn a coloring $C$ into a repeat multi-hypergraph $\mc{H}$. We must be careful for a few reasons. First, as we discussed earlier, a single coloring $C$ can give rise to many different $\mc{H}$. Second, although we saw that if $C$ is a $(p, q_{\text{quad}})$-coloring then $\mc{H}$ must be $(2\ell, \ell-1)$-free, in general the converse does not hold. In particular, if some set of vertices $S$ does not contain the edge $e_c$ then $S$ could have many repeats in the color $c$ but not span even one edge of $\mc{H}$. We will get around these issues by ensuring that our coloring uses each color at most twice, and we never use the same color on adjacent edges. For such a coloring $C$, the repeat multi-hypergraph $\mc{H}$ is unique. Furthermore, such a coloring $C$ is a $(p, q_{\text{quad}})$-coloring if and only if $\mc{H}$ is $(2\ell, \ell-1)$-free.

\begin{claim}\label{clm:quadupper}
\[
    f\rbrac{n, p, q_{\text{quad}}} \le \binom n2 - H^{(4)}(n; 2\ell, \ell-1).
\]
\end{claim}
\begin{proof}
We will construct a $\rbrac{p, q_{\text{quad}}}$-coloring that uses $\binom n2 - H^{(4)}(n; 2\ell, \ell-1)$ colors. 
   We start with an extremal graph $\mc{H}$ for the $H^{(4)}(n; 2\ell, \ell-1)$ problem. In other words, $\mc{H}$ has $n$ vertices, $H^{(4)}(n; 2\ell, \ell-1)$ edges, and properties \eqref{item:1}--\eqref{item:3}.

We construct a coloring as follows. Start with an edge $h_1\in \mc{H}$ and choose two arbitrary, disjoint pairs $e_{1}, f_{1} \subseteq h_1$ and assign them the color $c_{1}$. Assign all other pairs in $h_1$ a different color. Let the set of ``active'' pairs after step 1 be $A_1=\{e_1,f_1\}$.
Then define  
$$
H_1=\{h\in E(\mc{H})\setminus \{h_1\}: e  \subseteq h \text{ for some } e\in A_1\}.
$$ 
\begin{figure}[h]
\begin{tikzpicture}[scale=.8]
	\node (u) at (0,0) [vert,label=below:] {};
	\node (v) at (2,0) [vert,label=below:] {};
	\node (w) at (4,0) [vert,label=below:] {};
	\node (x) at (0,2) [vert,label=above:] {};
	\node (y) at (2,2) [vert,label=above:] {};
	\node (z) at (4,2) [vert,label=above:] {};
	\draw [dotted, thick] (u) -- (v);
	\draw [dotted, thick] (x) -- (y);
    \draw [very thick,red] (y) -- (z) node [midway,below,fill=white] {$e_2$};
    \draw [very thick,red] (v) -- (w) node [midway,above,fill=white] {$f_2$};
    \draw [blue, very thick] (x) -- (u) node [midway,right,fill=white] {$f_1$};
    \draw [dotted, thick ] (x) -- (v);
    \draw [dotted, thick] (y) -- (u);
    \draw [very thick, blue] (y) -- (v) node [midway,left=-4.7pt,] {$e_1$};
    \draw [dotted, thick] (y) -- (w);
    \draw [dotted, thick] (z) -- (v);
    \draw [dotted, thick] (z) -- (w);
    
    \draw[rounded corners] (-.5, -.5) rectangle (2.5, 2.5) {};
    \draw[rounded corners] (1.5, -.5) rectangle (4.5, 2.5) {};
    
    \node (h1) at (-.5,2.5) [label=above:$h_1$] {};
    \node (h2) at (4.5,2.5) [label=above:$h_2$] {};

\end{tikzpicture}
\caption{The coloring after step 2. Dotted edges are given all unique colors.}
\end{figure}
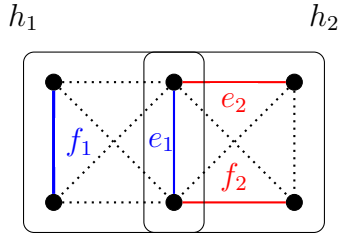
In general, assume that we have defined colors in $h_1, h_2, \ldots, h_{k-1}$ such that 
\begin{itemize}
\item 
$A_{k-1}= \{e_1, \ldots, e_{k-1}\} \cup \{f_1,\ldots,f_{k-1}\}$

\item 
$
H_{k-1}=\{h\in E(\mc{H})\setminus \{h_1, \ldots, h_{k-1}\}: e  \subseteq h \text{ for some } e\in A_{k-1}\}
$, and 
\item $|\bigcup_{i=1}^{k-1} h_i| = 2k$.
\end{itemize}
Notice that these are true for step 1. We choose an arbitrary edge $h_k$ from $H_{k-1}$. Notice that $\left(\bigcup_{i=1}^{k-1} h_i  \right) \cap h_k = e$ for some $e\in A_{k-1}$. Indeed, $e$ is clearly a subset of the expression and by property \eqref{item:2} if the cardinality of the intersection were $3$, then $\bigcup_{i=1}^{k} h_i$ has $2k+4-3 = 2k+1$ vertices that induces at least $k$ edges, violating the $(2k+1,k)$-free condition in $\mc{H}$. Thus, $\left|\bigcup_{i=1}^{k} h_i\right|=2(k+1)$. Then pick two disjoint, uncolored pairs $e_{k},f_{k}\subseteq h_k$ (there are two such choices) and color them $c_{k}$ and color all other uncolored pairs in $h_k$ a unique color. Finally, define 
$$
A_{k}= A_{k-1} \cup \{e_k,f_k\}
$$ 
and
$$
H_{k}=\{h\in E(\mc{H})\setminus \{h_1, \ldots, h_{k}\}: e  \subseteq h \text{ for some } e\in A_{k}\}.
$$
Continue in this way until $H_k = \emptyset$ for some $k$. Notice that $k < \ell-1$ since otherwise there would be a set of $2\ell$ vertices on $\ell-1$ edges, violating property \eqref{item:3}.  Then repeat this process with any edge other than $h_1, \ldots, h_k$. 
Continue until all edges have been processed, and give any uncolored pairs a unique color. 
Notice that each newly chosen edge will intersect the union of any of the former edges by at most $2$ by property \eqref{item:2}.
Note that for each edge in $\mc{H}$, the coloring has exactly one repeat, and there are $H^{(4)}(n,2\ell,\ell-1)$ edges. Thus when considering the coloring of the pairs, we obtain a coloring $C$ of $K_n$ with 
$$
|C| = \binom{n}{2} - H^{(4)}(n,2\ell,\ell-1).
$$
To verify that $C$ is a $(p, q_{\text{quad}})$-coloring, choose any set $S$ of $p=2\ell$ vertices. Clearly, in $\mc H$, the set $S$ induces at most $\ell-2$ hyperedges. And our coloring defines exactly one repeat per hyperedge, and none elsewhere. So the total number of distinct colors among the edges of $K_n[S]$ is at least $q_{\text{quad}} = \binom{2\ell}{2}-\ell+2$.
\end{proof}

Finally observe that Theorem \ref{thm:quadexact} follows from Lemma \ref{lem:shangguan} and Claims \ref{clm:FG}, \ref{clm:quadlower} and \ref{clm:quadupper}.

\subsection{Proof of \eqref{eqn:estquad} from Theorem \ref{thm:est}} In this subsection we establish the explicit bounds~\eqref{eqn:estquad}. They will follow from Theorem \ref{thm:quadexact} together with explicit bounds for the function $F^{(4)}$.

As we mentioned before, Delcourt and Postle \cite{DP22} proved some very general and powerful results to the effect that certain hypergraphs have almost-perfect matchings which avoid certain forbidden submatchings. Similar results were independently proved by Glock, Joos, Kim, K\"uhn and Lichev \cite{GJKKL}. Each team of researchers was motivated in part by finding approximate designs of high ``girth''. In particular, it follows just as well from either \cite{DP22} (Theorem 1.3) or \cite{GJKKL} (Theorem 1.1) that for any $\ell \ge 3$ there exists a linear 4-uniform hypergraph $\mc{H}$ on $n$ vertices with $\frac 1{12} n^2 +o(n^2)$ edges which is also $(2\ell, \ell-1)$-free. In other words, $F^{(4)}(n; 2\ell, \ell-1) \ge \frac 1{12} n^2 +o(n^2)$. Now the upper bound in \eqref{eqn:estquad} follows from Theorem~\ref{thm:quadexact}. 

We move on to the lower bound in \eqref{eqn:estquad}, which will follow from an upper bound on $H^{(4)}(n; 2\ell, \ell-1)$ (which of course gives an upper bound on $F^{(4)}(n; 2\ell, \ell-1)$ by~\eqref{eqn:FH}). In particular, we will be done when we prove the following:
\begin{claim} For $\ell \ge 2$ we have
\[
H^{(4)}(n; 2\ell, \ell-1) \le \frac{\ell - 2}{10\ell - 18} n^2 + o(n^2). 
\]  
\end{claim}
\noindent
The proof is a straightforward adaptation of Delcourt and Postle's proof of their Lemma 1.9 in \cite{DP22}.
\begin{proof}
    Let $\mc{H}$ be a 4-uniform hypergraph on $n$ vertices which is $(2\ell, \ell-1)$-free and $(2i+1, i)$-free for $i=2, \ldots, \ell-2$ (recall Definition \ref{def:H}). 
    
    Define a graph $G$ with $V(G)=E(\mc{H})$, where $e_1e_2 \in E(G)$ whenever $|e_1 \cap e_2| \ge 2$. Each component of $G$ must have order at most $\ell-2$ since $\mc{H}$ is $(2\ell, \ell-1)$-free. Let $\{e_1, \ldots e_b\}$ be a component in $G$ for some $1 \le b \le \ell-2$. Assume that the ordering $\{e_1, \ldots e_b\}$ is chosen so that for each $2\le i\le b$ there is $1\le j\le i-1$ such that $|e_i \cap e_j|\ge 2$.
    We claim that for each $i \ge 2$, $e_i$ has two vertices (in $V(\mc{H})$) which are not in $e_1 \cup \ldots \cup e_{i-1}$; otherwise, we would have a $(2i+1, i)$-configuration in $\mc{H}$. On the other hand, due to our choice of the ordering, there is an edge $e_j\in\{e_1 \cup \ldots \cup e_{i-1}\}$ such that $|e_i\cap e_j|\ge 2$. Consequently, $e_i$ shares only one pair of vertices with $e_1 \cup \ldots \cup e_{i-1}$ and so $e_i$ contains five pairs of vertices which are not subsets of any $e_{j}, j <i$. Of course $e_1$ contains six pairs and each edge after that has five more, so the total number of pairs contained in some $e_j, j \le b$ is at least $5b+1$.  Note that for two edges of $\mc{H}$, if they are in different components of $G$ then they do not share any pair of vertices in $\mc{H}$. 
    
    For $1 \le b \le \ell-2$ let $C_b$ be the number of components of $G$ of order $b$. Then we have
\begin{equation}\label{eqn:bcb}
     |E(\mc{H})| = \sum_{1 \le b \le \ell-2} b C_b,
\end{equation}
   which implies
   \begin{equation}\label{eqn:cb}
       \sum_{1 \le b \le \ell-2}  C_b \ge \frac{1}{\ell-2} |E(\mc{H})|.
   \end{equation}
   But now by summing the vertex-pairs in $\mc{H}$ we have
   \[
   \binom n2 \ge \sum_{1 \le b \le \ell-2} (5b + 1)  C_b \ge \rbrac{5 + \frac1{\ell-2}} |E(\mc{H})|,
   \]
where the last inequality uses \eqref{eqn:bcb} and \eqref{eqn:cb}. It follows that
\[
|E(\mc{H})| \le \frac{\ell - 2}{10\ell - 18} n^2 + o(n^2),
\]
which completes the proof. 
\end{proof}
\section{Linear threshold}\label{sec:lin}

In this section we address the linear threshold. First we prove Theorem \ref{thm:linlower} and the lower bound in \eqref{eqn:estlin}.

\subsection{Proof of Theorem \ref{thm:linlower} and the lower bound in \eqref{eqn:estlin}}

We start by comparing $F^{(3)}$ with $G^{(3)}$. This is analogous to Claim \ref{clm:FG}. 

\begin{claim}\label{clm:GFlin}
For all $p\ge 3$, 
$$
G^{(3)}(n; p, p-2) = F^{(3)}(n;p,p-2)+O(n).
$$
\end{claim}
\begin{proof}
Let $\mathcal H$ be an extremal multi-hypergraph for the $G^{(3)}(n;p,p-2)$ problem. So $\mathcal H$ has $G^{(3)}(n;p,p-2)$ edges and is $(p,p-2)$-free. 
Then $\mathcal H$ has at most $Cn$ edges of multiplicity at least 2, where $3C=2x$ and $x$ is whichever of $(p-2)/2$ or $(p-1)/2$ is an integer. 
Suppose to the contrary. Let $\mathcal H_2$ be the multi-hypergraph with $V(\mathcal H_2)=V(\mathcal H)$ and all edges from $\mathcal H$ with multiplicity at least 2. Then the average degree in $\mathcal H_2$ is at least $3C=2x$.
So there is a set of $2x$ edges on at most $1+2x$ vertices. But $2x \ge p-2$ and $2x+1\le p$, so this contradicts the fact that $\mathcal H$ is $(p,p-2)$-free. 

Now we form $\mathcal H'$ by deleting the edges with multiplicity at least $2$ that appear in $\mathcal H$. Since $\mathcal H'$ is also $(p,p-2)$-free, then it has at most $F^{(4)}(n;p,p-2)$ edges. In addition, we must delete at most $Cn$ edges of $\mathcal H$ to obtain $\mathcal H'$, so 
$$
F^{(3)}(n; p, p-2) \le G^{(3)}(n; p, p-2) \le F^{(3)}(n; p, p-2)+Cn,
$$
proving the claim. 
\end{proof}

The next claim is analogous to Claim \ref{clm:quadlower}.

\begin{claim}\label{clm:linlower}
For all $p \ge 3$
\begin{equation}
    f\rbrac{n, p, q_{\text{lin}}} \ge n -1- \frac 1n G^{(3)}(n; p, p-2).\nn
\end{equation}
\end{claim}
\begin{proof}
      Consider any $(p, q_{\text{lin}})$-coloring using color set $C$. So any set of $p$ vertices spans at most $p-3$ repeats.  We form the 3-uniform hypergraph $\mc{H}$ as follows. For each vertex $v$ and color $c$ used on at least one edge adjacent to $v$, say  $\{e_1, \ldots e_\ell\}$ is the set of edges adjacent to $v$ and colored $c$. Then $\mc{H}$ has the edges $e_1 \cup e_i$ for $i=2, \ldots, \ell$. 

      $\mc{H}$ is $(p, p-2)$-free, but it might have multiple edges which come from monochromatic triangles. Therefore
\[
|E(\mc{H})|  \le G^{(3)}(n; p, p-2).
\]
Each hyperedge of $\mc{H}$ corresponds to two edges of the same color sharing a vertex $v$, and so some particular vertex $v$ plays that role at most
\[
\frac 1n |E(\mc{H})|  \le \frac 1n G^{(3)}(n; p, p-2)
\]
times. But these hyperedges of $\mc{H}$ count all of the color repeats on edges incident with $v$. Thus the number of different colors used on edges adjacent with $v$ is at least 
\[
n-1 - \frac 1n G^{(3)}(n; p, p-2).
\]
\end{proof}

Theorem \ref{thm:linlower} now follows from Claims \ref{clm:GFlin} and \ref{clm:linlower}. In turn, the lower bound in \eqref{eqn:estlin} follows from Theorem \ref{thm:linlower} and Lemma 1.9 from Delcourt and Postle \cite{DP22}, which states that
\[
F^{(3)}(n, p, p-2) \le \frac{p-3}{4p-10}n^2 + o(n^2).
\]

\subsection{Proof of the upper bound in \eqref{eqn:estlin}}

We now turn to the upper bound at the linear threshold found in Theorem~\ref{thm:est}. We use the forbidden submatchings method of Delcourt and Postle~\cite{DP22b}. To introduce this method, we require the following definitions from~\cite{DP22b}.  

\begin{definition}
\begin{enumerate}[label=(\roman*)]
\ \\[-12pt]
\item The $i$-degree of a vertex $v$ of $H$, denoted $d_{H,i}(v)$, is the number of edges of $H$ of size $i$ containing $v$. The maximum $i$-degree of $H$, denoted $\Delta_i(H)$, is the maximum of $d_{H,i}(v)$ over all vertices $v$ of $H$. 
\item Let $G$ be a hypergraph. We say a hypergraph $H$ is a configuration hypergraph for $G$ if $V (H) = E(G)$ and $E(H)$ consists of a set of matchings of $G$ of size at least 2. We say a matching of G is $H$-avoiding if it spans no edge of $H$.
\item We say a hypergraph $G = (A, B)$ is \textit{bipartite with parts $A$ and $B$} if $V (G) = A \cup B$ and every edge of $G$ contains exactly one vertex from $A$. We say a matching of $G$ is \textit{$A$-perfect} if every vertex of $A$ is in an edge of the matching. We say a matching in $G$ is \textit{$H$-avoiding} if it contains no edge of $H$.
\item Let $H$ be a hypergraph. The \textit{maximum $(k,\ell)$-codegree} of $H$ is $$\Delta_{k,\ell}(H)=\max_{S\in \binom{V(H)}{\ell}} |\lbrace e \in E(H) : S \subset e, |e|=k\rbrace|$$
\end{enumerate}
\end{definition}

 We will use a slightly simplified version of Theorem 1.16 of Delcourt and Postle \cite{DP22b}. In particular the full version allows $H$ to have edges of size 2, at the cost of having to check some extra conditions. For our application we will avoid edges of size 2 in $H$.
 
\begin{theorem}[Delcourt and Postle \cite{DP22b}]\label{thm:dp}
For all integers $r\ge 2, g \ge 3$ and real $\beta \in  (0, 1)$, there exist an integer $D_\beta > 0$ and real $\alpha > 0$ such that following holds for all $D \ge D_\beta$:

Let $G = (A, B)$ be a bipartite $r$-bounded hypergraph  such that
\begin{enumerate}[label=(G\arabic*)]
\item every pair of vertices is in at most $D^{1-\beta}$ edges, and \label{thm-dp:A}
\item every vertex in $A$ has degree at least $(1+D^{-\alpha})D$ and every vertex in $B$ has degree at most $D$. \label{thm-dp:B}
\end{enumerate}
Let $H$ be a $g$-bounded configuration hypergraph of $G$ whose edges all have size at least 3 such that 
\begin{enumerate}[label=(H\arabic*)]
\item  $\Delta_i(H) \le \alpha \cdot D^{i-1} \log D$ for all $3\le i \le g$; \label{thm-dp:C}
\item $\Delta_{k,\ell}(H) \le D^{k-\ell-\beta}$ for all $3 \le \ell < k \le g$; and \label{thm-dp:D}
\end{enumerate}
Then there exists an $H$-avoiding $A$-perfect matching of $G$.
\end{theorem}

\begin{proof}[Proof of upper bound in \eqref{eqn:estlin}]
Fix some $\beta \in (0, 1)$, set $r=3$ and $g=\binom p2$ and let $\alpha>0$ be the value guaranteed by Theorem \ref{thm:dp}. Fix some $\eps$ with $0 < \eps <\alpha$. Let $C$ be a set of $n+n^{1-\eps}$ colors. 
Let $G=(A,B)$ be a bipartite hypergraph with parts $A=E(K_n)$ and $B=\left\lbrace v_c : c \in V(K_n), c\in C\right\rbrace$, and with edge set $$E(G)=\left\lbrace \lbrace uv, u_c, v_c\rbrace : u,v \in V(K_n), c\in C\right\rbrace.$$ Note that $G$ is $3$-uniform (and thus $3$-bounded). We intend to find an $A$-perfect matching $M$ in $G$, which will give us a $(p, p-2)$-coloring as follows. For each edge $\{uv, u_c, v_c\} \in M$ we just color the edge $uv$ with the color $c$. Since $M$ is $A$-perfect, every edge of $K_n$ gets exactly one color. Note that since $M$ is a matching, no two incident edges $uv$ and $vw$ in $K_n$ can get the same color $c$. 

We now define $H$, our configuration hypergraph of $G$. Of course we let $V(H)=E(G)$. Suppose $S\subseteq E(G)=V(H)$ is a matching, so $S$ corresponds to a coloring $c_S$ of some of the edges of $K_n$. We will let $S$ be an edge of $H$ if we have that
\begin{quote}
the number of vertices of $K_n$ spanned by edges that are colored by $c_S$ is $V(S)$ where $4 \le V(S) \le p$ and the number of color repetitions in $c_S$ is $R(S) \ge V(S)-2$.
\end{quote}
If any edge in $E(H)$ is not minimal (i.e.,  it properly contains another edge) we remove it. When $k=p$, an edge of $H$ corresponds to a violation of the $(p,q_{\text{lin}})$-condition in $K_n$. Including the edges of $H$ corresponding to $4 \le k \le p-1$ is important to verify the conditions Theorem~\ref{thm:dp}. It is easy to see that $H$ is $\binom{p}{2}$-bounded. Also, note that for all $e\in E(H)$, $|e|\ge 4$.

We now verify that condition~\ref{thm-dp:A} holds. Define $D = n$. Let $x,y \in V(G)$. Clearly, if $x,y \in A$, then the codegree is zero since there is exactly one member of $A$ in each edge of $G$. If $x \in A$ and $y \in B$ then $x=uv$ for some $u,v \in K_n$. If $y=u_c$ or $v_c$, then the codegree is $1$. Otherwise, the codegree is 0. Finally, if $x,y \in B$, then the codegree is either 0 or 1, depending on whether they share the same color subscript $c$. Thus, all codegrees in $G$ are at most $1$, verifying condition~\ref{thm-dp:A}. 

Next, we verify that condition~\ref{thm-dp:B} holds. For any vertex $uv\in A$, the degree of $uv$ in $G$ is exactly $|C|=n+n^{1-\eps}=D(1+D^{-\eps}) \ge D(1+D^{-\alpha})$. In addition, for any vertex $u_c\in B$, the degree of $u_c$ in $G$ is exactly $n-1 \le D$. So condition~\ref{thm-dp:B} is verified. 

Next, we verify that condition~\ref{thm-dp:C} holds. Let $e=\{uv, u_c, v_c\} \in V(H)=E(G)$. We count edges $I$ of $H$ of size $i$ with $e \in I$. For some $4 \le k \le p$, the number $V(I)$ of vertices of $K_n$ spanned by edges of $I$ in $G$ is $k$ and the number of color repetitions $R(I)$ is at least $k-2$. So besides $u$ and $v$, $V(I)$ must contain exactly $k-2$ other vertices of $K_n$. Let $x$ be the number of colors induced by $I$ other than $c$. We count $R(S)$ by taking the difference between the number of edges colored by $c_S$ and the number of distinct colors used. Thus we have $ i-(1+x) \ge k-2$, so $x\le i-k+1$. There are at most $\binom p2 = O(1)$ choices for remaining, unaccounted for, colored edges of $K_n$ edges in $I$. $I$ is determined by choosing $k-2$ vertices of $K_n$ and at most $i-k+1$ colors, so 
$$
\Delta_i(H) = O\left(\sum_{k=4}^p  n^{k-2} \cdot n^{i-k+1}\right) = O( n^{i-1}) \le \alpha \cdot D^{i-1} \log D
$$
for all $2\le i \le g$, verifying condition~\ref{thm-dp:C}. 
 
To verify condition~\ref{thm-dp:D}, fix $k$ and $\ell$, and let $L \subseteq V(H)$ have size $\ell$. We count the number of $K \in E(H)$ such that $L \subset K$ and $|K|=k$. If $V(L) > p$ there is no possible $K$, so we assume $V(L) \le p$. If $V(L)$ is 2 or 3 then $R(L)=0$. Otherwise we have $V(L) \ge 4$. If $R(L) \ge V(L)-2$ then there is no possible $K \subseteq L$ since we removed nonminimal edges from $H$, so assume $R(L) \le V(L)-3$. 
Let us count possible edges $K$ such that $V(K)=t$. Since $K$ is an edge of $H$, we have $R(K) = t-2$.  Then $V(L)-3 \ge R(L) = \ell - |C_S|$, so $|C_S| \ge \ell-s_1+3$. Suppose there are $\ell - R(L)$ many colors used by the coloring $c_L$, and say there are $x$ colors used by $c_K$  that are not used by $c_L$. Then the number of colors used by $c_K$ is 
\[
x + \ell - R(L) \le k - R(K) = k - t+2
\]
and so 
\[
x \le R(L) - \ell + k-t+2 \le V(L) - \ell + k - t - 1.
\]
To determine $K$ we choose $t-V(L)$ vertices of $K_n$ which are not touched by the coloring $c_L$, and then we choose $x$ many colors. Given that choice there are only a constant number of ways to choose which edges are colored and which colors they get. Therefore, 
$$
\Delta_{k,\ell}(H) \le O \left(\sum_{t \le p} n^{t-V(L)} \cdot n^{V(L) - \ell + k - t - 1}\right) = O(n^{k-\ell-1})\le D^{k-\ell-\beta}
$$
for all $2 \le \ell < k \le g$, verifying condition~\ref{thm-dp:D}.

Therefore, there exists an $H$-avoiding $A$-perfect matching of $G$, which corresponds to a $\left(p,\binom p2-p+3\right)$-coloring of $K_n$ using our set $C$ of $n+n^{1-\eps}$ colors. 
 
\end{proof}

\section{Concluding remarks}\label{sec:con}

In a previous draft of this paper we claimed that $f(n, 5,8) \ge \frac78 n + o(n)$ and conjectured that a matching upper bound holds. However there was an error in our proof, and recently in \cite{BCD24} we proved that instead we have $f(n, 5,8) = \frac67 n + o(n)$, matching the lower bound of G\'omez-Leos, Heath, Parker, Schwieder, and Zerbib \cite{GHPSZ}.

It is plausible to conjecture the following. Currently it is known only for $p=3, 4$. 

\begin{conjecture}
    The limit 
    \[
    \lim_{n \rightarrow \infty} \frac{f(n, p, q_{lin})}{n}
    \]
    exists for all $p\ge 3$.
\end{conjecture}
Finally, at the quadratic threshold, recall that we only proved that the analogous limit exists when $p$ is even. The same should likely also hold when $p$ is odd, as well as when $q$ is above the quadratic threshold.
\begin{conjecture}
    The limit 
    \[
    \lim_{n \rightarrow \infty} \frac{f(n, p, q)}{n^2}
    \]
    exists for all $p\ge 4$ and $q \ge q_{quad}$.
\end{conjecture}

\bibliographystyle{abbrv}

\end{document}